\newcommand{\C}{{\mathbb C}}           
\newcommand{\Z}{{\mathbb Z}}           
\newcommand{\N}{I\!\!N}                
\newcommand{\F}{{\mathbb F}}           
\newcommand{\R}{I\!\!R}                
\newcommand{\kpar}{K_{\rm par}}     
\newcommand{\Cpar}{\C_{\rm par}}
\newcommand{\lpar}{\lambda_{\rm p}} 

\documentclass[twoside,letterpaper,draft,11pt]{amsart}
\usepackage{amsmath}               
\usepackage{amsthm}                

\title{Partial Representations and Partial Group Algebras}
\author[M.\ Dokuchaev]{Michael Dokuchaev}
\address{Departamento de Matem\'atica, Universidade de S\~ao Paulo, Brazil}
\email{{\rm dokucha@ime.usp.br}}
\author[R.\ Exel]{Ruy Exel}
\address{Departamento de Matem\'atica, Universidade Federal de Santa Catarina, Brazil}
\email{{\rm exel@mtm.ufsc.br}}
\author[P.\ Piccione]{Paolo Piccione%
\thanks{Partially sponsored by CNPq, Brazil, Processo n.\ 301410/95-0 (RN)}}
\address{Departamento de Matem\'atica, Universidade de S\~ao Paulo, Brazil}
\email{{\rm piccione@ime.usp.br}}

\theoremstyle{definition}\newtheorem{defparrep}{Definition}[section]
\theoremstyle{definition}\newtheorem{defisoparrep}[defparrep]{Definition}
\theoremstyle{plain}\newtheorem{equivalence}[defparrep]{Proposition}
\theoremstyle{definition}\newtheorem{defkpar}[defparrep]{Definition}
\theoremstyle{remark}\newtheorem{remdim}[defparrep]{Remark}
\theoremstyle{plain}\newtheorem{representation}[defparrep]{Theorem}
\theoremstyle{plain}\newtheorem{isomorphic}[defparrep]{Corollary}

\theoremstyle{plain}\newtheorem{propab}{Proposition}[section]
\theoremstyle{plain}\newtheorem{structure}[propab]{Theorem} 
\theoremstyle{plain}\newtheorem{lattice}[propab]{Corollary}

\theoremstyle{plain}\newtheorem{prop}{Proposition}[section]
\theoremstyle{plain}\newtheorem{cor}[prop]{Corollary}
\theoremstyle{plain}\newtheorem{lemfloripa}[prop]{Lemma}
\theoremstyle{plain}\newtheorem{theo}[prop]{Theorem}
\theoremstyle{plain}\newtheorem{remintegraldomain}[prop]{Corollary}
\theoremstyle{remark}\newtheorem{counter}[prop]{Remark}

\theoremstyle{remark}\newtheorem{Zp}{Example}[section]
\theoremstyle{remark}\newtheorem{ZpZp}[Zp]{Example}
\theoremstyle{remark}\newtheorem{Zpq}[Zp]{Example}
\theoremstyle{remark}\newtheorem{Zpp}[Zp]{Example}
\theoremstyle{remark}\newtheorem{S3}[Zp]{Example}

\begin{document}

\begin{abstract}
The partial group algebra of a group $G$ over a field $K$, 
denoted by $\kpar(G)$, is the algebra whose representations  
correspond to the partial representations of $G$ over
$K$-vector spaces.
In this paper we study the structure of the partial
group algebra $\kpar(G)$, where $G$ is a finite group.
In particular, given two finite abelian groups $G_1$ and $G_2$,
we prove that if the characteristic of $K$ 
is zero, then $\kpar(G_1)$ is isomorphic to $\kpar(G_2)$
if and only if $G_1$ is isomorphic to $G_2$.
\end{abstract}

\maketitle


\begin{section}{Introduction}\label{sec:intro}
Let $G$ be a group, $K$ be a field and let $V$ be a $K$-vector space.  
By a partial representation of $G$ on $V$ we mean a map
  \[
  \pi : G \longmapsto {\rm End}(V)
 \]
  sending the unit group element to the identity operator on $V$, such
that for all $s,t\in G$ one has that
  $\pi(s)\pi(t)\pi(t^{-1})=\pi(st)\pi(t^{-1})$ and
  $\pi(s^{-1})\pi(s)\pi(t)=\pi(s^{-1})\pi(st)$.  Every representation
of $G$ on $V$ obviously fits this definition but there are several
interesting examples of partial representations which are not
representations in the usual sense of the word.

In order to study the representations of $G$ on $K$-vector spaces one
usually considers the group algebra $KG$, which is an associative
algebra whose representation theory is identical to the representation
theory of $G$.  Likewise, there is an associative algebra, called the
partial group algebra of $G$, and denoted $\kpar(G)$, which governs
the partial representations of $G$.  $\kpar(G)$ may be defined as
being the universal $K$-algebra generated by a set of symbols
$\left\{[t]\colon t\in G\right\}$, with relations
$[e]=1$, $[s^{-1}][s][t]=[s^{-1}][st]$, and
$[s][t][t^{-1}]=[st][t^{-1}]$, for every $t$ and $s$ in $G$.

Contrary to the case of $KG$, even if $G$ is an abelian group,
$\kpar(G)$ may be a non-commutative algebra.  For example, the complex
partial group algebra $\Cpar(\Z_4)$ of the cyclic group of order four is
isomorphic (see \cite{E1}) to
  \[7\,\C\oplus M_2(\C)\oplus M_3(\C),\]
  where by $7\,\C$ we mean the direct sum of seven copies of $\C$.
  On the other hand, $\Cpar(\Z_2\times\Z_2)$ is given by
  \[11\,\C\oplus M_3(\C).\]

The dimension of $\kpar(G)$ is  $2^{n-2}(n+1)$, where $n$ is
the order of $G$, and hence it is much larger than the dimension of  $KG$.
Partly for this reason $\kpar(G)$ keeps much more information about
the structure of $G$ than does $KG$.  In particular we have seen
that $\Cpar(\Z_4)\not\simeq \Cpar(\Z_2\oplus\Z_2)$ while the complex group algebra
of both groups of order four are isomorphic to $4\, \C$.

The fact that the partial group algebra construction is able to
differentiate between the two groups of order four, raises the
following question: under which circumstances
the isomorphism $\kpar(G) \simeq \kpar(H)$, implies
$G\simeq H$?  One of the main goals
of this paper is to prove that the above implication
holds provided that $G$ and $H$ are finite abelian groups and 
$K$ is a field (even integral domain) of characteristic zero; 
we also give a counterexample in the case of nonabelian groups. 

The history of the famous isomorphism problem for group rings
shows that, in general, $KG$ does {\em not\/} 
determine a group $G$ up to  isomorphism. Indeed, E.~C.~Dade in~\cite{Dade}
has constructed two nonisomorphic finite groups $G$ and $H$ such that $KG\simeq KH$
for all fields $K$. On the other hand, if $\Z G\simeq\Z H$, where $\Z$
is the ring of integers, then $KG\simeq KH$ for all fields $K$
(see~\cite[p.\ 664]{Passman}). So, the crucial question
is whether a finite group $G$ is determined up to isomorphism by
the group ring $\Z G$. 

Many deep positive results have been obtained on this problem, 
(see~\cite{Passman,RoggTay,Sandling,Sehgal}), however, a counterexample
has been recently announced by Martin Hertweck.
In view of this situation, it is natural to test the partial group rings for
the isomorphism question.

The concept of partial group representation was introduced in \cite{E1} and \cite{QR}, 
motivated by the
desire to study algebras generated by partial isometries on a Hilbert
space $H$.
  By a partial isometry on $H$ one usually means a bounded operator
$S\colon H\to H$ satisfying the equation $SS^*S=S$, where the star
refers to the Hilbert space adjoint.  Partial isometries clearly
include the isometries, i.e.~operators $S$ for which $S^*S=1$, as well
as unitary operators.

Partial isometries may be characterized geometrically as the linear
operators $S$ for which
  $\Vert S(\xi)\Vert = \Vert \xi\Vert$
  for all $\xi$ in the orthogonal complement of ${\rm Ker}(S)$.
It is therefore easy to see that the isometries are precisely the
injective partial isometries, while the unitary operators are those
which are bijective.

While the set of unitary operators forms a group under composition,
nothing like this holds for the set of all isometries, let alone the
case of partial isometries.  In particular the product of two
isometries is often not an isometry.  For this reason, algebraic
properties of partial isometries are extremely difficult to
investigate.  In particular, given a set of partial isometries, the
algebra they generate is often a wild object unlikely to yield to any
attempt at understanding its structure.

The first among the few well understood examples of algebras generated
by isometries is the Toeplitz algebra \cite{D1}, defined to be the
closed *-subalgebra of operators on the Hilbert space $l_2(\Z_+)$,
generated by the forward shift operator.  The Toeplitz algebra has
important applications in the study of Wiener--Hopf operators
\cite{R}, in $K$-theory \cite{PV}, and in $K$-homology \cite{BDF}.

Another thoroughly studied algebra generated by isometries is the so
called Cuntz algebra ${\mathcal O}_n$ \cite{C1} which is also relevant in
$K$-theory \cite{C2}.  ${\mathcal O}_n$ is a special case of the class of
Cuntz--Krieger algebras ${\mathcal O}_A$ defined as follows \cite{CK}:
given an $n\times n$ matrix $A = \{a_{ij}\}_{1\leq i,j\leq n}$ with
entries in $\{0,1\}$ one defines ${\mathcal O}_A$ as being the universal
$C^*$-algebra generated by partial isometries $S_1,\ldots,S_n$ subject
to the conditions:\smallskip

\begin{itemize}
\item[{CK$_1$)}] $\displaystyle\sum_{i=1}^n S_i S_i^* = 1$, and
  \smallskip 
\item[{CK$_2$)}] $\displaystyle S_i^* S_i = \sum_{j=1}^n a_{i,j} S_j S_j^*$.
\end{itemize}
  \smallskip

  The properties of ${\mathcal O}_A$ are closely related to the properties
of the Markov subshift whose transition matrix is given by $A$.

Attempting to study sets of isometries which generate a reasonably
tame algebra, Nica \cite{N}
studied isometric representations
  $
  t \mapsto V_t
  $
  of the positive cone $P$ of a  quasi-lattice ordered group $G$
satisfying a certain
``covariance condition''.  Nica's condition is satisfied in examples
related to the $C^*$-algebras generated by one-parameter groups of
isometries studied by Douglas \cite{D2}, as well as to the Cuntz
algebra, and hence it gives a satisfactory framework for the study of
algebras generated by isometries.  But it cannot deal with
Cuntz--Krieger algebras since the canonical set of generators of the
latter is formed by partial isometries rather than isometries.

The first indication that the concept of partial group representations
is relevant to the study of the Cuntz--Krieger algebras appeared in
\cite{E2}, where it is shown that, given any finite set
$S_1,\ldots,S_n$ of partial isometries on a Hilbert space $H$
satisfying CK$_1$--CK$_2$ above, there exists a partial representation
of the free group $\F_n$ sending the $i^{th}$ canonical generator to
$S_i$.  This idea was subsequently generalized in \cite{EL} to treat
the case of infinite matrices and was used to give the first
definition of  Cuntz--Krieger algebras for
transition matrices on infinitely many states, dropping the row-finite
condition imposed by Kumjian, Pask, Raeburn, and Renault in
\cite{KPRR}.

Also, it was proved by Quigg and Raeburn in \cite{QR} that every
isometric representation satisfying Nica's covariance condition extends to a
partial representation of the ambient group.
  Therefore there is a partial representation  underlying each
and every one of the examples of algebras generated by partial
isometries that have been successfully studied so far.

We should also note that by \cite{E3}, a necessary and sufficient
condition for a self-adjoint set $\mathcal S$ of partial isometries to be
closed under multiplication is that $\mathcal S$ be range of a partial
representation of a group.

\end{section}

\begin{section}{Basic definitions}\label{sec:basic}
Let $G$ be a finite group and $K$ be any field. We denote by $e$ the
identity of $G$; the notation $H\le G$ will mean that $H$ is a subgroup
of $G$; $H<G$ will mean that $H$ is a {\em proper\/} subgroup of $G$, i.e.,
$H\ne G$.  If $H\le G$, the symbol $(G:H)$ will denote the index of $H$ in $G$; 
by $[G,G]$ we will mean the commutator subgroup of  $G$, which is
the subgroup of $G$ generated by the elements of the form $sts^{-1}t^{-1}$, with
$s,t\in G$. For $H\le G$, ${\mathcal N}(H)$ will denote the {\em normalizer\/}
of $H$ in $H$, which is the maximal subgroup of $G$ containing $H$ as a normal
subgroup.
\begin{defparrep}\label{thm:defparrep}
A {\em partial representation\/} of $G$ on the $K$-vector space $V$ is a map 
$\pi:G\longmapsto{\rm End}(V)$ such that, for all $s,t\in G$, we have:
\begin{itemize}
\item[(a)] $\pi(s)\pi(t)\pi(t^{-1})=\pi(st)\pi(t^{-1})$;
\item[(b)] $\pi(s^{-1})\pi(s)\pi(t)=\pi(s^{-1})\pi(st)$;
\item[(c)] $\pi(e)=1$,
\end{itemize}
where $1={\rm Id}_V$ is the identity map on $V$. In general, a map
$\pi$ from $G$ to any unital algebra $A$ will be called a partial
representation of $G$ in $A$ if it satisfies (a), (b) and (c) above.
\end{defparrep}
\noindent
In other words, $\pi$ is a partial representation of $G$ if the equality
$\pi(s)\pi(t)=\pi(st)$ holds when the two sides are multiplied either by $\pi(s^{-1})$
on the left or by $\pi(t^{-1})$ on the right. In particular, every representation
of $G$ is a partial representation; moreover, if $H$ is any subgroup
of $G$ and $\pi:H\longmapsto{\rm End}(V)$ is a partial representation
of $H$, then the map $\tilde\pi:G\longmapsto{\rm End}(V)$ given by:
\[\tilde\pi(g)=\left\{\begin{array}{cl}\pi(g),&\text{if}\ g\in H;\\
0,&\text{otherwise}
\end{array}\right.\]
defines a partial representation of $G$. Obviously, Definition~\ref{thm:defparrep}
makes sense also for infinite groups.
\smallskip

The concept of equivalence for partial representations is the following.
Two partial representations $\pi_i:G\longmapsto {\rm End}(V_i)$, $i=1,2$, are
said to be {\em equivalent\/} if there exists a $K$-vector space  isomorphism
$\phi:V_1\longmapsto V_2$ such that $\phi\circ\pi_1(g)=\pi_2(g)\circ\phi$ for all 
$g\in G$. There is also a natural notion of {\em invariant subspace\/} for a
partial representation of a finite group $G$ and one can define the concept
of {\em irreducibility\/} for a partial representation in the obvious way.

A {\em unitary\/} representation of a group $G$ on an inner product space  $H$,
e.g., a Hilbert
space, is a representation $\pi$ of $G$ in the space $B(H)$
of all linear bounded operators on $H$ such that $\pi(g)$ is a unitary 
operator on $H$ for all $g\in G$.  
This means that $\pi(g)^*=\pi(g)^{-1}$ for all $g$, where
${}^*$ denotes the adjoint operation in $B(H)$.

The analogous concept for partial representations is the following:
\begin{defisoparrep}\label{thm:defisoparrep}
A {\em partial isometric representation\/} of $G$ on the complex Hilbert
space $H$ is a map $\pi:G\longmapsto B(H)$ such that:
\begin{itemize}
\item[(a)] $\pi(s)\pi(t)\pi(t^{-1})=\pi(st)\pi(t^{-1})$;
\item[(b)] $\pi(t^{-1})=\pi(t)^*$;
\item[(c)] $\pi(e)=1$;
\end{itemize}
for all $s,t\in G$.
\end{defisoparrep}
It is not hard to see that, if $\pi$ is a partial isometric representation
of $G$ on $H$, then $\pi(g)$ is a partial isometry in $B(H)$, i.e.,
$\pi(g)\pi(g)^*$ is an orthogonal projection. 

In~\cite{E1}, Definition~\ref{thm:defisoparrep} is taken as the definition
of a partial representation of a group $G$. We will show in a moment that the two
definitions coincide, up to equivalence.

In first place, observe that if $\pi:G\longmapsto B(H)$ is a partial isometric
representation, then $\pi$ is a partial representation. Namely, using properties
(a) and (b) of Definition~\ref{thm:defisoparrep}, we get:
\[\pi(s^{-1})\pi(s)\pi(t)=(\pi(t^{-1})\pi(s^{-1})\pi(s))^*=(\pi(t^{-1}s^{-1})\pi(s))^*=
\pi(s^{-1})\pi(st),\]
for all $s,t\in G$. Conversely, we have the following:
\begin{equivalence}\label{thm:equivalence}
Let $K=\R$ or $\C$.
Given any partial representation $\pi$ of $G$ on a $K$-vector space,
there exists a partial isometric representation $\pi_1$ of $G$
which is equivalent to $\pi$.
\end{equivalence}
\begin{proof}
We will prove that, given any partial representation $\pi:G\longmapsto{\rm End}(V)$,
with $V$ a complex vector space, then there exists a positive
definite $K$-valued inner product
$\langle\cdot,\cdot\rangle$ in $V$ such that:
\begin{equation}\label{eq:equivunit}
\langle\pi(g)\xi,\eta\rangle=\langle \xi,\pi(g^{-1})\eta\rangle,
\end{equation}
for all $g\in G$ and all $\xi,\eta\in V$. In order to prove this, we argue as follows.

For all $t\in G$, let $\epsilon(t)$ denote the element $\pi(t)\pi(t^{-1})$; 
it is  easy to check that
the $\epsilon(t)$'s form a family of commuting idempotents in ${\rm End}(V)$:
\[\epsilon(t)^2=\pi(t)\,\pi(t^{-1})\,\pi(t)\,\pi(t^{-1})=\pi(t)\,\pi(e)\,\pi(t^{-1})=
\epsilon(t);\]
and, for all $t,s\in G$,
\begin{equation}\label{eq:2a}
\begin{split}
\pi(t)\,\epsilon(s)&=\pi(t)\,\pi(s)\,\pi(s^{-1})=\pi(ts)\,\pi(s^{-1})=
\pi(ts)\,\pi(s^{-1}t^{-1})\,\pi(ts)\,\pi(s^{-1})\!=\\
&=\pi(ts)\,\pi(s^{-1}t^{-1})\,\pi(t)=\epsilon(ts)\,\pi(t),
\end{split}
\end{equation}
from which we compute:
\begin{equation}\label{eq:1b}
\epsilon(t)\,\epsilon(s)=\pi(t)\,\pi(t^{-1})\,\epsilon(s)=\!\pi(t)\,\epsilon(t^{-1}s)\,\pi(t^{-1})=
\epsilon(s)\,\pi(t)\,\pi(t^{-1})=\!\epsilon(s)\,\epsilon(t).
\end{equation}
Given a subset $S$ of $G$, we denote by $P_S$ the element in ${\rm End}(V)$
given by:
\begin{equation}\label{eq:PS}
P_S=\prod_{s\in S}\epsilon(s)\prod_{s\not\in S}(1-\epsilon(s)).
\end{equation}
Using (\ref{eq:2a}), it is easily seen
that $\pi(t)P_S=P_{tS}\,\pi(t)$. Observe that, since
$\epsilon(1)=1$, 
$P_S$ is zero unless $S$ contains the identity $e$ of $G$; moreover, if $t\not\in S$,
then $P_S\,\pi(t)=0$, because $\epsilon(t)\pi(t)=\pi(t)$.  
Also, the following identity holds:
\begin{equation}\label{eq:sommaPS}
\sum_{S\subseteq G}P_S=1,
\end{equation}
where the sum is taken over all possible subsets of $G$, including the empty subset. 
In order to prove (\ref{eq:sommaPS}), it suffices to observe that, since the $\epsilon$'s
commute, we have the following combinatorial formula:
\[
1=\prod_{s\in G}1=\prod_{s\in G}\big(1-\epsilon(s)+\epsilon(s)\big)=
\sum_{S\subseteq G}\left(\Big[\prod_{t\in S}\epsilon(t)\Big]\cdot\Big[\prod_{t\not\in
S}(1-\epsilon(t))\Big]\right).
\]

Now, let $[\cdot,\cdot]$ be any $K$-valued positive definite 
inner product in $V$; we consider the new inner product
$\langle\cdot,\cdot\rangle$ on $V$, given by:
\begin{equation}\label{eq:new}
\langle\xi,\eta\rangle=\sum_S\sum_{t\in G}\big[P_S\,\pi(t)\xi, P_S\,\pi(t)\eta\big].
\end{equation}
In order to prove that $\langle\cdot,\cdot\rangle$ is positive definite, 
since $\big[\cdot,\cdot\big]$ is positive definite, 
it suffices to show that $P_S\xi=0$ for all $S\subseteq G$ implies $\xi=0$.
This follows easily from~(\ref{eq:sommaPS}).

We now prove  (\ref{eq:equivunit}). To see this, we start by observing that, for
all $S\subseteq G$, $t,r\in G$ and $\xi,\eta\in V$, the following equality holds:
\begin{equation}\label{eq:PSpi1}
\big[P_S\,\pi(t)\pi(r)\xi,P_S\,\pi(t)\eta\big]=\big[P_S\,\pi(tr)\xi, P_S\,\pi(t)\eta\big].
\end{equation}
Namely, observe that if $t\not\in S$, then $P_S\,\pi(t)=0$ and both sides of
(\ref{eq:PSpi1}) vanish. On the other hand, if $t\in S$, then $P_S\epsilon(t)=
P_S$, hence:
\[P_S\,\pi(t)\pi(r)=P_S\,\epsilon(t)\pi(t)\pi(r)=P_S\,\epsilon(t)\pi(tr)=
P_S\,\pi(tr),\]
which proves (\ref{eq:PSpi1}). Similarly, by symmetry one checks that:
\begin{equation}\label{eq:PSpi2}
\big[P_S\,\pi(t)\xi,P_S\,\pi(t)\pi(r)\eta\big]=\big[P_S\,\pi(t)\xi,P_S\,\pi(tr)\eta\big].
\end{equation}
For $g\in G$ and $\xi,\eta\in V$, using (\ref{eq:PSpi1}) and (\ref{eq:PSpi2})
we now compute as follows:
\begin{equation*}
\begin{split}
\langle\pi(g)\xi,\eta\rangle&=\sum_{S\subseteq G}\sum_{t\in G}\big[P_S\,\pi(t)\pi(g)\xi,
P_S\,\pi(t)\eta\big]=\\
&=\sum_{S\subseteq G}\sum_{t\in G}\big[P_S\,\pi(tg)\xi,
P_S\,\pi(t)\eta\big]=\quad\mbox{(setting $t'=tg$)}\\
&=\sum_{S\subseteq G}\sum_{t'\in G}\big[P_S\,\pi(t')\xi,
P_S\,\pi(t'g^{-1})\eta\big]=\\
&=\sum_{S\subseteq G}\sum_{t'\in G}\big[P_S\,\pi(t')\xi,
P_S\,\pi(t')\pi(g^{-1})\eta\big]=\langle\xi,\pi(g^{-1})\eta\rangle,
\end{split}
\end{equation*}
which concludes the proof.
\end{proof}

The group algebra $KG$ is the algebra with the same representation theory
of the group $G$; in a similar fashion, we can define the {\em partial group
algebra\/} $\kpar(G)$, whose representations are in one-to-one correspondence
with the partial representations of $G$ as follows.
\begin{defkpar}\label{thm:defkpar}
Given a finite group $G$ and a field $K$,
the partial group $K$-algebra $\kpar(G)$ is the universal $K$-algebra with unit $1$
generated by the set of symbols $\{[g]:g\in G\}$, with relations:
\begin{itemize}
\item[(1)] $[e]=1$;
\item[(2)] $[s^{-1}][s][t]=[s^{-1}][st]$;
\item[(3)] $[s][t][t^{-1}]=[st][t^{-1}]$;
\end{itemize}
for all $s,t\in G$.
\end{defkpar}
It is straightforward to check that, if ${\mathcal A}$ is any unital
$K$-algebra and $\pi:G\longmapsto {\mathcal A}$ is a partial representation of $G$ in
${\mathcal A}$, then $\pi$ extends uniquely by linearity to a representation $\phi:\kpar(G)
\longmapsto {\mathcal A}$ such that $\phi([g])=\pi(g)$. 

Conversely, if $\phi:\kpar(G)\longmapsto
A$ is a unital homomorphism of $K$-algebras, then $\pi(t)=\phi([t])$ gives a partial
representation of $G$ in ${\mathcal A}$. 

Rather than giving generators and relations, we want to give an alternative
(and constructive) description of the algebra $\kpar(G)$.

To this aim, associated to $G$, we will consider a finite groupoid 
(Brandt groupoid)
$\Gamma=\Gamma(G)$, 
whose elements are pairs $(A,g)$, where $g\in G$ and $A$ is a  subset 
of $G$ containing the identity $e$ and the element $g^{-1}$. 
We will refer to \cite{Re} for the basic definitions and properties of groupoids. 
Observe that $e,g\in gA$.
\smallskip

The multiplication of pairs  $(A,g)\cdot(B,h)$
in $\Gamma$ is defined  for pairs for which   $A=hB$, in which
case we set:
$$(hB,g)\cdot(B,h)=(B,gh);$$
the inverse of the pair $(A,g)$ is $(gA,g^{-1})$.  
The {\em source\/} and the {\em range\/}  map on 
$\Gamma$ are respectively $s(A,g)=(A,e)$ and  $r(A,g)=(gA,e)$.
The set of {\em units\/} of $\Gamma$ is denoted by $\Gamma^{(0)}$, the 
$\vert\cdot\vert$ function will denote the {\em cardinality\/} of a set.
\medskip

We recall that if $\Gamma$ is any groupoid, denoting by $\Gamma^{(2)}
\subseteq\Gamma\times\Gamma$ the set of its composable pairs,
the {\em groupoid algebra\/} $K\Gamma$ is a $K$-vector space
linearly generated by the elements  
of $\Gamma$, and with the multiplication given by:
\[\gamma_1\cdot\gamma_2=\left\{
\begin{array}{ll}\gamma_1\gamma_2,&\text{if}\ (\gamma_1,\gamma_2)\in\Gamma^{(2)}, 
\\ 0,  &  \text{otherwise}, 
\end{array}
\right.\]
 and extended linearly on $K\Gamma$.

Let $K\Gamma(G)$ denote the $K$-algebra of the groupoid $\Gamma(G)$;  
the dimension of this algebra is equal to the cardinality of $\Gamma(G)$,
which is the number of pairs $(A,g)$ as described above. If $\vert G\vert=n$,
it is easily computed:
\begin{equation}\label{eq:dimension}
{\rm dim}(K\Gamma(G))=\sum_{k=0}^{n-1}(k+1)\binom{n-1}{k}=2^{n-2}(n+1). 
\end{equation}
\begin{remdim}\label{thm:remdim}
Observe that the right hand side of (\ref{eq:dimension}) is a strictly 
increasing function on $n$. In particular, if $G$ and $H$ are finite groups
such that $K\Gamma(G)$ is isomorphic to $K\Gamma(H)$, then $\vert G\vert=\vert H\vert$.
\end{remdim}

The elements of the form $(A,e)$ are {\em idempotents\/} in $K\Gamma(G)$, they are
mutually orthogonal and their sum is the identity of $K\Gamma(G)$:
\[(A,e)^2=(A,e),\qquad(A,e)(B,e)=0\quad\hbox{if $A\ne B$},\qquad\sum_{A\ni e}(A,e)=1. \]

We define the  map $\lpar:G\longmapsto K\Gamma(G)$ by:
\begin{equation}\label{eq:lpar}
\lpar(g)=\sum_{A\ni g^{-1}}(A,g).
\end{equation} 
Observe that in the above formula we are summing over all possible subsets
$A$ of $G$ that contain the elements $e$ and $g^{-1}$.

We claim that $\lpar$ is a partial representation of $G$. Namely, for $g,h\in G$,
we compute:
\begin{eqnarray*}
\lpar(g^{-1})\lpar(g)\lpar(h)&=&
\sum_{\substack{ g\in A \\ g^{-1}\in B \\ h^{-1}\in C}}
(A,g^{-1})(B,g)(C,h)=\\&=&
\sum_{\substack{h^{-1}\in C\\ g^{-1}\in hC}}
(ghC,g^{-1})(hC,g)(C,h)=\\
&=&\sum_{\substack{h^{-1}\in C\\ h^{-1}g^{-1}\in C}}(C,h).
\end{eqnarray*}
On the other hand, we have:
\begin{eqnarray*}
\lpar(g^{-1})\lpar(gh)&=&\sum_{\substack{g\in A\\ h^{-1}g^{-1}\in C}}
(A,g^{-1})(C,gh)=\\
&=&\sum_{\substack{h^{-1}\in C\\ h^{-1}g^{-1}\in C}}
(C,h).
\end{eqnarray*}
Hence, $\lpar(g^{-1})\lpar(g)\lpar(h)=\lpar(g^{-1})\lpar(gh)$ for all $g,h\in G$.
Similarly, one proves the equality $\lpar(g)\lpar(h)\lpar(h^{-1})=\lpar(gh)\lpar(h^{-1})$.
Moreover, by definition
\[\lpar(e)=\sum_{A\ni e}(A,e)=1,\]
and the claim is proven.
We are ready to show that the partial representation theory of $G$ is the same as
the representation theory of the algebra $K\Gamma(G)$:
\begin{representation}\label{thm:representation}
There is a one-to-one correspondence between the partial representations
of $G$ and the representations of $K\Gamma(G)$. More precisely, if ${\mathcal A}$ is any unital
$K$-algebra, then $\pi:G\longmapsto {\mathcal A}$ is a partial representation of $G$ if and
only if there exists a unital algebra homomorphism $\tilde\pi:K\Gamma(G)\longmapsto {\mathcal
A}$ such that $\pi=\tilde\pi\circ\lpar$. Moreover, such homomorphism $\tilde\pi$ is
unique. 
\end{representation}
\begin{proof}
If $\tilde\pi:K\Gamma(G)\longmapsto {\mathcal A}$ is a unital homomorphism of
$K$-algebras, then clearly $\pi=\tilde\pi\circ\lpar:G\longmapsto {\mathcal A}$
is a partial representation of $G$ in ${\mathcal A}$.

Conversely, suppose that $\pi:G\longmapsto {\mathcal A}$ is a partial representation
of $G$. Arguing as in the proof of Proposition~\ref{thm:equivalence}, 
for all $r\in G$ we denote by
$\epsilon(r)$ the element of ${\mathcal A}$ given by $\pi(r)\pi(r^{-1})$. 
Recall from (\ref{eq:2a}) and (\ref{eq:1b}) that the $\epsilon(r)$'s are commuting 
idempotents and
\begin{equation*}
\pi(s)\epsilon(r)=
\epsilon(sr)\pi(s)
\end{equation*}
for all $s,r\in G$. 
Similarly
\[\epsilon(r)\pi(s)=\pi(s)\epsilon(s^{-1}r).\]
For all elements $(A,g)$ in $\Gamma(G)$, we define:
\begin{equation}\label{eq:defpitilde}
\tilde\pi(A,g)=\pi(g)\left(\prod_{r\in A}\epsilon(r)\right)\left(
\prod_{s\not\in A}(1-\epsilon(s))\right).
\end{equation}
For $(A,g)$ and $(B,h)$ in $\Gamma(G)$, we have:
\begin{equation*}
\begin{split}
\tilde\pi(A,g)&\,\tilde\pi(B,h)=\\& =\pi(g)
\prod_{r\in A}\epsilon(r)  
\cdot\prod_{s\not\in A}(1-\epsilon(s))\cdot\pi(h)\cdot
\prod_{t\in B}\epsilon(t)\cdot\prod_{v\not\in B}(1-\epsilon(v))=\\
& =\pi(g)\,\pi(h)\prod_{r\in A}\epsilon(h^{-1}r)  
\prod_{s\not\in A}(1-\epsilon(h^{-1}s))
\prod_{t\in B}\epsilon(t)\prod_{v\not\in B}(1-\epsilon(v))=\\
& =\pi(g)\,\pi(h)\prod_{r\in h^{-1}A}\epsilon(r)  
\prod_{s\not\in h^{-1}A}(1-\epsilon(s))
\prod_{t\in B}\epsilon(t)\prod_{v\not\in B}(1-\epsilon(v)).
\end{split}
\end{equation*} 
Now, if $h^{-1}A\ne B$, i.e.\ if $A\ne hB$, then, either there exists 
one element $r$ that belongs to $h^{-1}A$ but does not belong to $B$,
or there is an element $r$ that belongs to $B$ but not to $h^{-1}A$.
In both cases, the product $\tilde\pi(A,g)\tilde\pi(B,h)$ contains the factor 
$\epsilon(r)(1-\epsilon(r))=0$ and therefore it is null.

If $h^{-1}A=B$, then, since $h^{-1}\in h^{-1}A$, 
the product $\tilde\pi(A,g)\tilde\pi(B,h)$ becomes:
\begin{equation*}
\begin{split}
\tilde\pi(A,g)\,\tilde\pi(B,h)&= 
\pi(g)\,\pi(h)\prod_{r\in h^{-1}A}\epsilon(r)  
\prod_{s\not\in h^{-1}A}(1-\epsilon(s))=\\&= 
\pi(g)\,\pi(h)\,\pi(h^{-1})\,\pi(h)\prod_{%
\substack{r\in h^{-1}A\\ r\ne h^{-1}}}
\epsilon(r)  
\prod_{s\not\in h^{-1}A}(1-\epsilon(s))=\\&= 
\pi(gh)\prod_{r\in h^{-1}A}\epsilon(r)  
\prod_{s\not\in h^{-1}A}(1-\epsilon(s))=\tilde\pi((A,g)\cdot(B,h)).
\end{split}
\end{equation*}
This proves that $\tilde\pi$ is multiplicative on $\Gamma(G)$, and thus its
linear extension to  $K\Gamma(G)$ gives a homomorphism  of 
$K\Gamma(G)$ in $A$. 

We now prove that $\tilde\pi$ is unital. By definition, 
since $\pi(e)=1$ and recalling (\ref{eq:PS}) and (\ref{eq:sommaPS}), we have:
\begin{equation}\label{eq:311}
\begin{split}
\tilde\pi(1)&=\tilde\pi\left(\sum_{A\ni e}(A,e)\right)=\sum_{A\ni
e}\tilde\pi(A,e)= 
\sum_{A\ni e}\prod_{r\in A}\epsilon(r)\prod_{s\not\in A}(1-\epsilon(s))=\\
&=\sum_{A\ni e}P_A=\sum_{A\subseteq G}P_A=1.
\end{split}
\end{equation}

Furthermore, we have: 
\begin{eqnarray*}
\tilde\pi\circ\lpar(g)&=&\tilde\pi\left(\sum_{A\ni g^{-1}}
(A,g)\right)=\sum_{A\ni g^{-1}}\tilde\pi(A,g)=\\&=&\pi(g)\sum_{A\ni g^{-1}} 
\prod_{r\in A}\epsilon(r)\prod_{s\not\in A}(1-\epsilon(s))=\\&=&
\pi(g)\pi(g^{-1})\pi(g)\sum_{A\ni g^{-1}}\prod_{r\in A\setminus\{g^{-1}\}}
\epsilon(r)\prod_{s\not\in A}(1-\epsilon(s))=\\&=&
\pi(g)\sum_{A\ni g^{-1}} (\epsilon(g^{-1})+1-\epsilon(g^{-1}))
\!\prod_{r\in A\setminus\{g^{-1}\}}\!\epsilon(r)\prod_{s\not\in A}
(1-\epsilon(s))\!=\\&=&\pi(g)\sum_{B}\prod_{r\in B}\epsilon(r)\prod_{s\not\in B}   
(1-\epsilon(s))=\\&=&\pi(g)\tilde\pi(\sum_B(B,e))=\pi(g)\tilde\pi(1)=\pi(g).
\end{eqnarray*}
To conclude, we are left with the proof of the uniqueness of the
homomorphism $\tilde\pi$ satisfying $\tilde\pi\circ\lpar=\pi$.
In order to show this, it suffices to prove that the set $\lpar(G)$ 
generates the whole algebra $K\Gamma(G)$. 

To this aim, let $(B,h)$ be an arbitrary element of $\Gamma(G)$,
where \[B=\{b_1^{-1},b_2^{-1},\dots,b_{k-1}^{-1},h^{-1}\}\] is
a subset of $G$ containing the identity. The set of such pairs
$(B,h)$ form a vector space basis of $K\Gamma(G)$. Let's denote by
${\mathcal A}$ the subalgebra of $K\Gamma(G)$ generated by $\lpar(G)$.
Let $\{g_1,g_2,\dots,g_k\}$ be the family of elements of $G$ defined by:
\begin{eqnarray*}
&&g_1=b_1,\ g_2g_1=b_2,\ g_3g_2g_1=b_3,\ \cdots,\\
&&\ g_{k-1}g_{k-2}\cdot\dots
\cdot g_1=b_{k-1},\ g_kg_{k-1}\cdot\dots\cdot g_1=h.
\end{eqnarray*}
We consider the element in $\mathcal A$ given by the product
$\lpar(g_k)\lpar(g_{k-1})\cdots\lpar(g_1)$:
\[\lpar(g_k)\lpar(g_{k-1})\cdots\lpar(g_1)=
\sum_{\substack{A_1\ni g_1^{-1}\\ A_2\ni g_2^{-1}\\ \vdots\\
A_k\ni g_k^{-1}}}
(A_k,g_k)\cdot(A_{k-1},g_{k-1})\cdot\dots\cdot(A_1,g_1).\]
Using the rules of multiplication in $\Gamma(G)$, it is easily seen that
the sum equals:
\[\sum_{\substack{A_1\ni g_1^{-1}\\ g_1A_1\ni g_2^{-1}\\ \vdots
\\ g_{k-1}g_{k-2}\cdot\dots\cdot g_1A_1\ni g_k^{-1}}}
(A_1,g_kg_{k-1}\cdot\dots\cdot g_1)=\sum_{A_1\supseteq B}(A_1,h).\]
Hence, for every $(B,h)\in\Gamma(G)$, $\mathcal A$ contains the 
element:
\[\sum_{A\supseteq B}(A,h).\]
Now, suppose that $G\setminus B=\{x_1,x_2,\dots,x_N\}$.
We have:
\[\sum_{A\supseteq B}(A,h)-\sum_{A\supseteq B\cup\{x_1\}}(A,h)= 
\sum_{\substack{A\supseteq B\\ A\not\ni x_1}}(A,h),\]
which means that, for every $(B,h)\in\Gamma(G)$ and every $x\not\in 
B$,  $\mathcal A$ contains all elements of the form:
\[\sum_{\substack{A\supseteq B\\ A\not\ni x}}(A,h).\]
Now, 
\[\sum_{\substack{A\supseteq B\\ A\not\ni x_1}}
(A,h)-\sum_{\substack{A\supseteq B\cup\{x_2\}
\\ A\not\ni x_1}}(A,h)=\sum_{\substack{A\supseteq B\\ A\not\ni x_1,x_2
}} 
(A,h),\]
hence, for every $x_1,x_2\not\in B$, $\mathcal A$ contains the element 
\[\sum_{\substack{A\supseteq B\\ A\not\ni x_1,x_2}}(A,h). \]
We repeat the argument, using the fact that:
\[\sum_{\substack{A\supseteq B\\ A\not\ni x_1,x_2}}(A,h)
-\sum_{\substack{A\supseteq B\cup\{x_3\}\\ A\not\ni x_1,x_2}}(A,h)=
\sum_{\substack{A\supseteq B\\ A\not\ni x_1,x_2,x_3}}(A,h).\]
By induction, we see that $\mathcal A$ contains all elements of the form:
\[\sum_{\substack{A\supseteq B\\ A\not\ni x_1,x_2,\dots,x_N}}(A,h)=(B,h), \]
which proves that ${\mathcal A}=K\Gamma(G)$ and we are done.
\end{proof}
Note that the bijection between the partial representations of $G$
and the representations of $\kpar(G)$ given by the composition
with the map $\lpar$ preserves the notions of equivalence and irreducibility.
\begin{isomorphic}\label{thm:isomorphic}
The groupoid algebra $K\Gamma(G)$ is isomorphic to the partial group
algebra $\kpar(G)$.
\end{isomorphic}
\begin{proof}
The maps $[\;]:G\longmapsto\kpar(G)$, given by $g\longmapsto[g]$, and $\lpar:G\longmapsto
K\Gamma(G)$ are partial representations of $G$. 

By the universal properties of the two algebras
$\kpar(G)$ and $K\Gamma(G)$, there exists unital $K$-algebra homomorphisms
$\tilde\pi:K\Gamma(G)\longmapsto\kpar(G)$ and $\phi:\kpar(G)\longmapsto K\Gamma(G)$
such that $\tilde\pi(\lpar(g))=[g]$ and $\phi([g])=\lpar(g)$. 

It is easily checked that the composite maps $\tilde\pi\circ\phi$ and $\phi\circ\tilde\pi$
are the identities on their domains, as they are the identity on a set of generators,
hence $\phi$ and $\tilde\pi$ are isomorphisms.
\end{proof}
\end{section}

\begin{section}{On the structure of the partial group algebra}
\label{sec:structure}
In this section we study the structure of the groupoid algebra $K\Gamma(G)$;
we prove that $K\Gamma(G)$ is the direct sum of matrix algebras over the rings
$KH$, where $H$ is any subgroup of $G$. To this aim, we introduce the
following notation.

Given a finite group $H$ and a positive integer $m$, let $\Gamma_m^H$ denote the
groupoid whose elements are triplets 
$(h,i,j)$, where $h\in H$ and $i,j\in\{1,2,\dots,m\}$. The source
and the range map on $\Gamma_m^H$ are defined by $s(h,i,j)=j$ and
$r(h,i,j)=i$. The product is defined by $(g,i,j)\cdot(h,j,k)=(gh,i,k)$.
The units of $\Gamma_m^H$ are the elements of the form $(e,i,i)$, $i=1,\ldots,m$.

It is sometimes useful to represent a groupoid $\Gamma$ with an oriented graph
$E_\Gamma$, whose vertices are the units of the groupoids, and each element
$\gamma\in\Gamma$ gives an oriented edge of $E_\Gamma$ from the vertex
$s(\gamma)$ to the vertex $r(\gamma)$. A connected component of $E_\Gamma$ gives
a subgroupoid of $\Gamma$.

In the case of the groupoid $\Gamma_m^H$, the corresponding graph $E_{\Gamma_m^H}$
has precisely $m$ vertices, and between any two vertices there are precisely
$\vert H\vert$ oriented edges (in each direction) labeled by the elements
of $H$.

We have the following simple result:
\begin{propab}\label{thm:propab}
Let $\Gamma$ be a groupoid such that $E_\Gamma$ is connected
and $m=\vert\Gamma^{(0)}\vert$ is finite. Let $x_1$ be any vertex
of $E_\Gamma$, and  $H$ be the {\em isotropy
group\/} of $x_1$, which is defined by:
\[H=\{\gamma\in\Gamma:s(\gamma)=r(\gamma)=x_1\}.\]
Then,
\begin{enumerate}
\item[(a)] $\Gamma\simeq\Gamma_m^H$;
\item[(b)] $K\Gamma\simeq M_m(KH)$.
\end{enumerate}
\end{propab}
\begin{proof}
Let $x_1,\ldots,x_m$ be the units of $\Gamma$; for all $i,j=1,\ldots,m$ define
${\mathcal E}_{i,j}$ by:
\[{\mathcal E}_{i,j}=\{\gamma\in\Gamma:s(\gamma)=x_i,\ r(\gamma)=x_j\}.\]
Since $\Gamma$ is connected, then ${\mathcal E}_{i,j}$ is non empty for all 
$i,j$. 
Namely, if $y_1=x_i,\ y_2,\ldots, y_{k+1}=x_j$ is a sequence
in $\Gamma^{(0)}$ and $\gamma_1,\ldots,\gamma_k$ is a sequence in $\Gamma$
such that $s(\gamma_i)=y_i$ and $r(\gamma_i)=y_{i+1}$ for all $i$,
then the element $\gamma=\gamma_k\cdot\gamma_{k-1}\cdot\ldots\cdot\gamma_1$
is in ${\mathcal E}_{i,j}$.


Clearly, $\Gamma$ is the (disjoint) union of all the ${\mathcal E}_{i,j}$'s.
Fix a family of elements
$\gamma_i\in {\mathcal E}_{1,i}$, with $i=1,\ldots,m$.  
Then, for every element $\gamma\in\Gamma$, $\gamma_j^{-1}\gamma\gamma_i=h\in H$
and $\gamma$ can be written {\em uniquely\/} in the form $\gamma=\gamma_j h \gamma_i^{-1}$ 
with $h\in H$, where $s(\gamma)=x_i$ and $r(\gamma)=x_j$. 
On the other hand, every element of the form $\gamma_j h \gamma_i^{-1}$,
with $h\in H$, belongs to ${\mathcal E}_{i,j}$. Observe that, in particular,
$\vert {\mathcal E}_{i,j}\vert$ is constant  equal to $\vert H\vert$.

The map $\gamma=\gamma_j h \gamma_i^{-1}\longmapsto (h,j,i)$ gives the desired 
isomorphism between $\Gamma$ and $\Gamma_m^H$, proving part (a).

Finally, let $(e_{i,j})_{i,j=1,\ldots,m}$ denote the set of  
matrix units of $M_m(K)$, i.e., $e_{i,j}$ is the matrix with
null entries, except for the entry in the position $(i,j)$ which is 
equal to the identity of $K$.
It is easy to check that, for any group $H$, the map
$(h,i,j)\longmapsto h\otimes e_{i,j}$ extends by linearity to a $K$-algebra
isomorphism between $K\Gamma_m^H$ and $M_m(KH)\simeq M_m(K)\otimes_K KH$,
which gives part (b) and concludes the proof.
%
%
\end{proof}

We remark that in Proposition~\ref{thm:propab} the group $H$
can be infinite.

If a groupoid $\Gamma$ is a (finite) disjoint union $\bigcup_i\Gamma_i$ of 
subgroupoids $\Gamma_i$, then the groupoid algebra $K\Gamma$ is the
direct sum $\bigoplus_iK\Gamma_i$; in particular, the groupoid algebra
$K\Gamma$ is a direct sum of algebras of the form $M_m(KH)$ where $m\ge0$ and
$H$ is a finite group. 

If $G$ is a finite group and $\Gamma(G)$ is the groupoid constructed in
Section~\ref{sec:basic}, for $1\le k\le\vert G\vert$, by the $k$-th {\em level\/}  
of $\Gamma(G)$ we will mean the subgroupoid of $\Gamma(G)$ consisting of pairs $(A,g)$ 
with $\vert A\vert=k$. Similarly, the $k$-th level of the graph $E_{\Gamma(G)}$
is meant to be the subgraph of $E_{\Gamma(G)}$ consisting of all vertices
$(A,e)$ and all edges $(A,g)$ with $\vert A\vert=k$. 

\smallskip

In the next result we give a more precise description of the algebra
$K\Gamma(G)$, and an algorithm for its computation. 

\begin{structure}\label{thm:structure}
The groupoid algebra $K\Gamma(G)$ is of the form:
\begin{equation}\label{eq:decomposition}
K\Gamma(G)=\bigoplus_{\substack{H\le G \\ \\ 1\le m\le(G:H)}}
c_m(H)\,M_m(KH),
\end{equation}
where $c_m(H)\,M_m(KH)$ means the direct sum of $c_m(H)$ copies
of $M_m(KH)$. Moreover, the following
recursive formula holds:
\begin{equation}\label{eq:recursive}
c_m(H)=\frac1m\cdot\Bigg(\binom{(G:H)-1}{m-1}-m\sum_{%
\substack{H<B\le G\\ \\ (B:H)\, \big\vert\,
m}}\frac{c_{\frac{m}{(B:H)}}(B)}{(B:H)}\Bigg).
\end{equation}
\end{structure}
\begin{proof}
Let $A$ be any subset of $G$ containing the identity of $G$. 
In the course of the proof, we will identify $A$ with the vertex $(A,e)$
of the graph $E_{\Gamma(G)}$.
We denote by $H=S(A)$
the {\em stabilizer\/} of $A$  in $G$, given by:
\[S(A)=\{g\in G:gA=A\}.\]
In the graph $E_{\Gamma(G)}$, $S(A)$ is identified with the set of
edges departing and terminating at the vertex $(A,e)$.
Observe that, since $e\in A$, then $H\subseteq A$. 

Since $H$ acts on the left on $A$, then $A$ is the union of right cosets
of $H$, say:
\[A=\bigcup_{i=1}^mHt_i,\qquad t_1=e,\]
where
\[m=\frac{\vert A\vert}{\vert H\vert}.\]
We now prove that the sub-groupoid of $\Gamma(G)$ corresponding
to the connected component of the vertex $A$ of the graph
$E_{\Gamma(G)}$ is (isomorphic to) the groupoid $\Gamma^H_m$.
First of all, we prove that the connected component of $A$ in $E_{\Gamma(G)}$
contains precisely $m$ vertices, given by the subsets of $G$ of the
form $A_i=t_i^{-1}A$, $i=1,\dots,m$. These sets are all distinct, in fact,
if $t_i^{-1}A=t_j^{-1}A$, then $t_jt_i^{-1}\in H$ and $Ht_i=Ht_j$,
which implies that  $i=j$. Moreover, if $r\in Ht_i$ for some $i$, say $r=ht_i$
for $h\in H$, then
$r^{-1}A=t_i^{-1}h^{-1}A=t_i^{-1}A$, which implies that the connected component
of $A$ contains no other vertex but the $A_i$'s.
Observe that the stabilizer of $A$, which is $H$, 
coincides with the isotropy group of the unit $(A,e)\in\Gamma(G)^{(0)}$.
Hence, by
Proposition~\ref{thm:propab} the algebra $M_m(KH)\simeq K\Gamma_m^H$ is a direct
summand of $K\Gamma(G)$ and  $K\Gamma(G)$ is the direct
sum of algebras arising from this construction.
This proves (\ref{eq:decomposition}).\smallskip

For the proof of the
second part of the thesis, we fix $H$ and 
$m$ as in the hypothesis. From the first part of the proof, it is clear that
in order to determine the number of occurrences of the algebra $M_m(KH)$ in
(\ref{eq:decomposition}) we need to count the number of vertices at the level 
$k=m\times\vert H\vert$  of $E_{\Gamma(G)}$ whose stabilizer is $H$. 

We define $b_m(H)$ to be the number of such vertices; moreover, 
we consider the index of $H$ in its normalizer ${\mathcal N}(H)$:
\[\delta(H)=({\mathcal N}(H):H)=\text{number of subgroups of $G$ conjugate to $H$}.\]

Since the conjugation is
an automorphism of $G$, if $H$ and $H'$ are conjugate in $G$ then, by symmetry,
$b_m(H)=b_m(H')$ for all $m$. Observe that two vertices in the
same connected component of $E_{\Gamma(G)}$ have conjugate stabilizers; hence,
the number of vertices at the level $k$ whose stabilizers
are conjugate to $H$ is precisely $\delta(H)\cdot b_m(H)$. 
The connected component of each of these vertices, that gives a contribution
of one copy of $M_m(KH)$ as a direct summand of $\kpar(G)$, contains $m$ vertices.
Hence, the number of summands  $M_m(KH)$ in $\kpar(G)$ that
arise from this construction is $\frac1m\cdot{b_m(H)\cdot\delta(H)}$. Finally,
by symmetry, the contribution of the single subgroup $H$ is:
\begin{equation}\label{eq:bmcm}
c_m(H)=\frac{b_m(H)\cdot\delta(H)}{m\cdot\delta(H)}=\frac{b_m(H)}m.
\end{equation}

Every  vertex at the $k$-th
level of $E_{\Gamma(G)}$ which is fixed by $H$ is a subset of $G$ (containing $e$)
which is union of $m$ right cosets of $H$ (including $H$); the number of such vertices
is $\binom{(G:H)-1}{m-1}$. Some of them will have stabilizer which is bigger
than $H$, and this number is clearly given by     
\[\sum_{\substack{H< B\le G\\ \\(B:H)\,\vert\,m}}
b_{\scriptscriptstyle{\frac{m}{(B:H)}}}(B),\]
hence we have the following recursive formula for the coefficients
$b_m(H)$:
\begin{equation}\label{eq:recursbm}
b_m(H)= \binom{(G:H)-1}{m-1}- \sum_{\substack{H< B\le G\\ \\(B:H)\,\vert\,m}}
b_{\scriptscriptstyle{\frac{m}{(B:H)}}}(B).
\end{equation}
From (\ref{eq:bmcm}) and (\ref{eq:recursbm}) we obtain immediately (\ref{eq:recursive}),
which concludes the proof.
\end{proof}
In the following corollary we give a sufficient condition for two finite groups
to have isomorphic partial group algebras:
\begin{lattice}\label{thm:lattice}
Let $G_1$ and $G_2$ two finite groups having isomorphic subgroup lattices,
such that the corresponding subgroups have isomorphic group rings over $K$.
Then, $\kpar(G_1)$ is isomorphic to $\kpar(G_2)$.
\end{lattice}
\begin{proof}
Under our hypotheses, there exists a bijection between the family of subgroups of 
$G_1$ and $G_2$ that preserves the order of the subgroups and their inclusions. 
Hence, \eqref{eq:recursive} implies that the coefficients $c_m(H)$ are also
preserved.  The conclusion follows from \eqref{eq:decomposition}, considering
that the corresponding subgroups have isomorphic group rings.
\end{proof}
The corollary suggests the following problem: given two finite groups
$G_1$ and $G_2$ with isomorphic partial group algebras over $K$, does there
exist an isomorphism of their subgroup lattices, such that the corresponding subgroups
have isomorphic group algebras over $K$?

\end{section}


\begin{section}{The Isomorphism Problem. Abelian Groups}
\label{sec:abelian}
As it was mentioned in the Introduction, in the case of the nonisomorphic groups
$\Z_2\times\Z_2$ and $\Z_{4}$, the corresponding partial group algebras 
are also nonisomorphic. Motivated by this
example,  in this section we prove that given any two finite abelian groups
$G_1$ and $G_2$, their partial group algebras $\kpar(G_1)$ and $\kpar(G_2)$
are isomorphic if and only if $G_1$ and $G_2$ are, provided that 
the characteristic of $K$ is zero.

We will show with a counterexample that for non commutative groups, in
general, the isomorphism problem for the partial group algebras has a 
negative answer.

Given a finite group $G$ and a positive integer $m\le\vert G\vert$, we define
$\beta_m(G)$ by:
\begin{equation}\label{eq:defbeta}
\beta_m(G)=\left\{\begin{array}{cl}
0,&\!\!\mbox{if $G$ does not have subgroups of order $m$;}\\ \\
\displaystyle\sum_{\substack{ H\le G\\ \vert H\vert=m }}(H:[H,H]),&
\mbox{otherwise.}
\end{array}\right.
\end{equation}
Moreover, for all $m\in\N$, we denote by $\gamma_m(G)$ the number of subgroups
of order $m$ of $G$.

Clearly, if all the subgroups of order $m$ of $G$ are abelian, then
$\beta_m(G)=m\gamma_m(G)$.
\begin{prop}\label{thm:prop}
Let $K$ be an algebraically closed field whose characteristic is zero.
Let $k\in\N$ be a divisor of $\vert G\vert$ with $\vert G\vert\ne2k$ and such that every prime
that divides $\vert G\vert$ also divides $\vert G\vert\cdot k^{-1}$. Then, the multiplicity
of the summand $M_{\vert G\vert k^{-1}-1}(K)$ in the Wedderburn decomposition of
$\kpar(G)$ is:
\begin{equation}\label{eq:multiplicity}
\frac k{\vert G\vert -k}
\sum_{\substack{m\,\big\vert\,\vert G\vert\\1\le m<\frac{k\vert G\vert}{\vert G\vert-k} 
}}
\binom{\frac{\vert G\vert}m-1}{\frac{\vert G\vert}k-2}\beta_m(G).
\end{equation}
\end{prop}
\begin{proof}
By Corollary~\ref{thm:isomorphic} and Proposition~\ref{thm:structure}, we have
an isomorphism:
\begin{equation}\label{eq:*}
\kpar(G)\simeq\bigoplus M_{n_i}(KH_i),
\end{equation}
where the $H_i$'s are (not necessarily distinct) subgroups of $G$. 

If $n_i<\frac{\vert G\vert}k-1$, then we claim that the summand $M_{n_i}(KH_i)$ does not
contain  a Wedderburn component isomorphic to $M_{\frac{\vert G\vert}k-1}(K)$. Indeed,
if it did, then the algebra $KH_i$ would contain a direct summand of the form
$M_s(K)$, where \[n_i\,s=\frac{\vert G\vert}k-1,\] 
and $s$ is the dimension of some irreducible $K$-representation of $H_i$.

Let $p>1$ be a prime dividing $s$;
in particular, $p$ divides $\frac{\vert G\vert}k-1$.
Then, since the characteristic of $K$ is zero, $s$ divides $\vert H_i\vert$ (see
e.g.~\cite[pag.\ 216]{CR}), and consequently it also  divides
$\vert G\vert$. By hypothesis it follows that $p$ must divide $\frac{\vert G\vert}k$. This
yields
$p=1$, which is a contradiction and proves our claim.

\noindent
Thus, the summands of (\ref{eq:*}) whose Wedderburn decomposition contains
$M_{\frac{\vert G\vert}k-1}(K)$ are precisely the ones given by $M_{n_i}(KH_i)$, with
$n_i=\frac{\vert G\vert}k-1$.

Observe that it must be:
\[n_i\,\vert H_i\vert<\vert G\vert;\]
for, the equality
\[\left(\frac{\vert G\vert}k-1\right)\,\vert H_i\vert=\vert G\vert\]
 would imply $\vert H_i\vert=\vert G\vert$ and $\frac{\vert G\vert}k-1=1$, which
would give $\vert G\vert=2k$, a contradiction with our hypothesis.

Now, each summand $M_{n_i}(KH_i)$ in (\ref{eq:*}), with $n_i=\frac{\vert G\vert}k-1$, 
contains exactly $( H_i:[H_i,H_i])$ Wedderburn components isomorphic to
$M_{\frac{\vert G\vert}k-1}(K)$ (see for instance \cite{S}). 
Therefore, in order to compute the multiplicity of the summand $M_{\frac{\vert G\vert}k-1}(K)$
in the Wedderburn decomposition of $\kpar(G)$, it suffices to determine the
multiplicity of $M_{\frac{\vert G\vert}k-1}(KH)$ in (\ref{eq:*}) for a {\em fixed\/}
subgroup $H$.

To this aim, let $H$ be a fixed subgroup of $G$ and
let $A$ be a subset of $G$ containing the identity, with
$\vert A\vert=(\frac{\vert G\vert}k-1)\vert H\vert$, and with stabilizer
$S(A)\supseteq H$.

If $S(A)\ne H$, take a prime $p$ that divides $(S(A):H)$. Since $\vert S(A)\vert$
divides $\vert A\vert$, then $p$ must divide $\frac{\vert G\vert}k-1$. But $p$ divides
$\vert G\vert$, hence $\frac{\vert G\vert}k$  by our hypothesis. Again, this implies
that $p=1$, which shows that $S(A)=H$. 

Consequently, every set $A$ which is the union of $\frac{\vert G\vert}k-1$ distinct
cosets of $H$ (including $H$) has stabilizer $H$. The total number of such sets is
easily computed as $\binom{(G:H)-1}{\frac{\vert G\vert}k-2}$,
so the multiplicity of $M_{\frac{\vert G\vert}k-1}(KH)$ in (\ref{eq:*}) with fixed $H$
is:
\[\frac k{\vert G\vert-k}\binom{(G:H)-1}{\frac{\vert
G\vert}k-2}.\]  It follows that the number of algebras $M_{\frac{\vert G\vert}k-1}(K)$
contributed into the Wedderburn decomposition of $\kpar(G)$ by \[\bigoplus_{\substack{%
H\le G\\ \vert H\vert=m}}M_{\frac{\vert G\vert}k-1}(KH) \]  
is given by:
\[\frac k{\vert G\vert-k}\binom{\frac{\vert G\vert}m-1}{\frac{\vert G\vert}k-2}\cdot
\beta_m(G).\]
Since the inequality $\left(\frac{\vert G\vert}k-1\right)\,m<\vert G\vert$ is equivalent
to $m<\frac{k\,\vert G\vert}{\vert G\vert-k}$, the result follows.
\end{proof}

\begin{cor}\label{thm:cor}
Let $G$ and $H$ be finite abelian groups  
with $\vert G\vert=\vert H\vert=p^n\cdot a>4$ ($p\, \vert^{\!\!\!-} a$)
for some prime $p$ and some integer $a$. Suppose that $K$ is an algebraically
closed field whose characteristic is zero such that
$\kpar(G)\simeq \kpar(H)$. Then, $\gamma_{p^j}(G)=\gamma_{p^j}(H)$
for all $j\in\N$ such that $2j\le n$.
\end{cor}
\begin{proof}
We prove more in general that if $k\in\N$ is a divisor of $\vert G\,\vert$
such that the following two conditions are satisfied:
\begin{itemize}
\item[(1)] every prime that divides $\vert G\vert$ divides also $\frac{\vert G\vert}k$;
\item[(2)] either $\vert G\vert>k^2$ or $G$ is a $p$-group and $\vert G\vert\ge k^2$,
\end{itemize}
then $\gamma_m(G)=\gamma_m(H)$ for all $m\le k$.
Under our hypotheses, the number $k=p^j$ satisfies (1) and (2) above.
For the sake of shortness, in the course of this proof, 
an integer number $k$ satisfying (1) will be called a {\em small\/} 
divisor of $\vert G\,\vert$.

We proceed by induction on $k$; the case $k=1$ is clearly trivial.

Let $k>1$ be fixed.

Observe first that the inequality $m<\frac{\vert
G\vert\,k}  {\vert G\vert-k}$ holds if and only if $m\le k$. Indeed,
\[\frac{\vert G\vert\,k}{\vert G\vert-k}=k+\frac{k^2}{\vert G\vert-k},\]
so $m\le k$ obviously implies $m<\frac{\vert G\vert\,k}{\vert G\vert-k}$. 
Now, if $G$ is a $p$-group and $\vert G\vert\ge k^2$, then $\vert G\vert\ge2k$, and
consequently $\frac{k}{\vert G\vert-k}\le1$. Thus, 
\[m<\frac{\vert G\vert\,k}{\vert G\vert-k}=k\,\left(1+\frac k{\vert G\vert-k}\right)
\le2k.\]
Since both $m$ and $k$ are powers of $p$, it follows that $m\le k$.

Now, if $G$ is any group with $\vert G\vert>k^2$, then $\frac{\vert G\vert}k\ge k+1$,
which gives $\frac{k^2}{\vert G\vert-k}\le1$. It follows
\[m<\frac{\vert G\vert\,k}{\vert G\vert-k}=k+\frac{k^2}{\vert G\vert-k}\le k+1,\]
and $m\le k$, as claimed.\smallskip

Let now $m\le k$ be fixed; observe that, by condition  (2) above,
since $\vert G\vert>4$ we have $\vert G\vert>2k$. Thus, 
using Proposition~\ref{thm:prop}, we have that the multiplicity of
$M_{\frac{\vert G\,\vert}k-1}(K)$ in the Wedderburn decomposition of $\kpar(G)$ is:
\begin{equation}\label{eq:**}
\frac k{\vert G\,\vert-k}\sum_{\substack{1\le m\le k\\ m\,\big\vert\,
\vert G\vert}}m\,\binom{\frac{\vert G\vert}m-1}{\frac{\vert G\vert}k-2}\,\gamma_m(G).
\end{equation}
By the uniqueness of the Wedderburn decomposition, the number in (\ref{eq:**}) is
invariant by $K$-algebra isomorphisms. Since the coefficient of $\gamma_k(G)$
in the formula (\ref{eq:**}) is non zero, to conclude the proof it suffices
to show that $\gamma_m(G)=\gamma_m(H)$ for all $m<k$.

If $m<k$ is a {\em small\/} divisor of $\vert G\,\vert$, then,
by our induction hypothesis, it follows  $\gamma_m(G)=\gamma_m(H)$.

Suppose now that $m<k$ is not a {\em small\/} divisor of $\vert G\,\vert$, and let
$p_1,p_2,\ldots,p_s$ be all the prime divisors of
$\vert G\vert$ that do {\em not\/} divide $\frac{\vert G\vert}m$. 
We prove that in this case there exists $\tilde m<k$ a {\em small\/}
divisor of $\vert G\,\vert$ such that 
$\gamma_m(G)=\gamma_{\tilde m}(G)$.

Since $G$ is abelian,
every subgroup of $G$ of order $m$ is a product of the form $P_1\times P_2\times\dots\times
P_s\times Q$, where $P_i$ is the (unique) Sylow $p_i$-subgroup of $G$ and $Q$ is some subgroup
of
$G$ such that none of the $p_i$'s divide $\vert Q\vert$. Set $\tilde m=\vert Q\,\vert$.

It is easily seen that $\gamma_m(G)=\gamma_{\tilde m}(G)$ and that $\tilde m$
is a {\em small\/} divisor of $\vert G\,\vert$. 

It follows that $\gamma_m(G)=\gamma_{\tilde m}(G)=\gamma_{\tilde
m}(H)=\gamma_m(H)$, and we are done. 
\end{proof}
We present one more preliminary result:
\begin{lemfloripa}\label{thm:lemfloripa}
Let  $P$ and $Q$ be abelian $p$-groups, $p$ prime, with $\vert P\vert=\vert Q\vert
=p^n$ for some integer $n$. Suppose that $\gamma_{p^j}(P)=\gamma_{p^j}(Q)$
for all $j$ such that $2j\le n$. Then, $P\simeq Q$.
\end{lemfloripa}
\begin{proof}
Let $P=\langle g_1\rangle\times\langle g_2\rangle\times\cdots\times\langle g_t\rangle$
and $Q=\langle h_1\rangle\times\langle h_2\rangle>\times\ldots\times\langle h_r\rangle$
be the cyclic factor decompositions of $P$ and $Q$, and set:
\[\Omega_k(P)=\big\{g\in P:g^{p^k}=1\big\}.\]
Clearly, $\gamma_k(P)=\gamma_k(\Omega_k(P))$ for all $k\ge0$.

We can obviously assume $n\ge2$, since the case $n=1$ is trivial.
The equality $\gamma_p(P)=\gamma_p(Q)$ implies that $\Omega_1(P)\simeq\Omega_1(Q)$ and,
in particular, $t=r$.

Now, suppose by contradiction that $P\not\simeq Q$. For sufficiently large $k$ we have
$\Omega_k(P)=P$ and $\Omega_k(Q)=Q$, so we can choose the minimal $j$
such that $\Omega_j(P)\not\simeq\Omega_j(Q)$. Clearly, $j\ge2$, because $\Omega_1(P)\simeq
\Omega_1(Q)$.

Let $n_k$ (respectively, $m_k$) the number of the $g_i$'s (respectively, of the $h_i$'s)
whose order is equal to $p^k$. Since $\Omega_j(P)\not\simeq\Omega_j(Q)$,
the numbers $n_{j-1}$ and $m_{j-1}$ must be different, say $n_{j-1}>m_{j-1}$.
Moreover, by the minimality of $j$, it must be $n_i=m_i$ for all
$i\le j-2$, because $\Omega_{j-1}(P)\simeq\Omega_{j-1}(Q)$. Therefore,
$\Omega_j(Q)$ can be obtained from $\Omega_j(P)$ by replacing some direct
cyclic factors of order $p^{j-1}$ with cyclic factors of order $p^j$.
But, obviously, such a replacement increases the total
number of subgroups of order $p^j$, from which it follows
$\gamma_j(\Omega_j(P))<\gamma_j(\Omega_j(Q))$, and, consequently,
$\gamma_j(P)<\gamma_j(Q)$.
Then, by our hypothesis, we must have $2j>n$.

On the other hand, $Q$ has at least one element $h_i$ whose order is greater or
equal to $p^j$, because $n_{j-1}<m_{j-1}$ and $t=r$. If $Q$ has only one such
$h_i$, then $n_{j-1}-m_{j-1}=1$, $n_j=n_{j+1}=\ldots=0$, and therefore
$Q$ can be obtained from $P$ by replacing a direct cyclic factor
of order $p^{j-1}$ with a cyclic factor of order greater or equal to $p^j$.
But this means that $\vert P\vert<\vert Q\vert$, which is impossible.
Hence, $Q$ contains at least two elements $h_i$ with order greater or equal
to $p^j$, which implies that $p^n=\vert Q\vert\ge p^{2j}$, and $2j\le n$,
a contradiction.
\end{proof}
We are now ready to prove our main result:
\begin{theo}\label{thm:theo}
Let $G$ and $H$ be two finite abelian groups and $K$ be a field
with characteristic  zero  such that
the partial group algebras $\kpar(G)$ and $\kpar(H)$ are isomorphic. 
Then, $G$ and $H$ are isomorphic groups.
\end{theo}
\begin{proof}
We can assume that $K$ is algebraically closed. Namely, if $L$ is any algebraically
closed field that contains $K$, then the isomorphism $\kpar(G)\simeq\kpar(H)$
implies the isomorphism \[{L_{\rm par}}(G)\simeq L\otimes_K\kpar(G)
\simeq L\otimes_K\kpar(H)\simeq {L_{\rm par}}(H).\]

As we observed in Remark~\ref{thm:remdim}, if $\kpar(G)\simeq\kpar(H)$, then
$G$ and $H$ have the same cardinality.

The case $\vert G\,\vert=4$ is treated separately.
The partial group algebra $\kpar(\Z_4)$ of the cyclic group of order $4$
and $\kpar(\Z_2\times\Z_2)$ of the Klein $4$-group are isomorphic  to
$7\,K\oplus M_2(K)\oplus M_3(K)$ and $11\,K\oplus M_3(K)$ respectively
(see Section~\ref{sec:examples}). Hence, we have $\kpar(\Z_4)\not\simeq\kpar(\Z_2\times
\Z_2)$.
These are the only (abelian) groups of order $4$, 
and so we may now assume that $\vert G\vert>4$.

From Corollary~\ref{thm:cor}, it follows that $\gamma_{p^j}(G)=\gamma_{p^j}(H)$ 
for all prime number $p$ and all integer $j$ such that $p^{2j}$ divides
$\vert G\,\vert$.
Let $G=P_1\times P_2\times\dots\times P_t$ where $P_i$ is the Sylow $p_i$-subgroup
of $G$. Then $H=Q_1\times Q_2\times\dots\times Q_t$, where $Q_i$ is the Sylow
$p_i$-subgroup of $H$; clearly $\vert
Q_i\vert=\vert P_i\vert$ for all $i=1,\ldots,t$.

Fix a $p_i=p$ and set $P=P_i$, $Q=Q_i$. Since $\gamma_{p^j}(G)=\gamma_{p^j}(P)$
and $\gamma_{p^j}(H)=\gamma_{p^j}(Q)$, Lemma~\ref{thm:lemfloripa} implies
that $P\simeq Q$. 

As $p=p_i$ is arbitrary, we finally obtain $G\simeq H$.
\end{proof}
The notion of partial $K$-representation and partial group
$K$-algebra can be obviously extended to the case of an integral domain $K$.
Since the proofs of Theorem~\ref{thm:representation}, Corollary~\ref{thm:isomorphic},
Proposition~\ref{thm:propab} and Theorem~\ref{thm:structure} do not involve
inverses of the coefficients in $K$, these facts remain true for integral
domains. Moreover, we have the following:
\begin{remintegraldomain}
\label{thm:remintegraldomain}
Let $G$ and $H$ be two finite abelian groups and $K$ be an integral domain
with characteristic zero. If
the partial group algebras $\kpar(G)$ and $\kpar(H)$ are isomorphic,  
then  $G$ and $H$ are isomorphic groups.
\end{remintegraldomain}
\begin{proof}
It suffices to observe  that if $\kpar(G)\simeq\kpar(H)$, 
then $L_{\rm par}(G)\simeq L_{\rm par}(H)$, where $L$ is 
the field of fractions of $K$. Hence,
the conclusion follows from Theorem~\ref{thm:theo}.
\end{proof}
\begin{counter}\label{thm:counter} Theorem~\ref{thm:theo} does {\em not\/}
hold for non commutative groups.
To see this, we present the following example.

Let $G_1$ and $G_2$ be the  
groups:
\[\begin{split}
&G_1=\langle a,b,c\;\big\vert\; a^{11}=b^{11}=c^5=1,\ ab=ba,\ c^{-1}ac=a^3,\
c^{-1}bc=b^9\rangle\\ &G_2=\langle a,b,c\;\big\vert\; a^{11}=b^{11}=c^5=1,\ ab=ba,\
c^{-1}ac=a^3,\ c^{-1}bc=b^4\rangle.
\end{split}\]
Then, $G_1$ and $G_2$ are nonisomorphic groups of order $605$
with isomorphic partial group algebras over an algebraically closed field $K$
of characteristic zero.
Indeed, according to a result of~\cite{Rott}, these two groups have isomorphic
subgroup lattices, and the corresponding subgroups have the same orders.
It is easily seen that, in fact, the corresponding proper subgroups are isomorphic
and, in particular, they have isomorphic group algebras. Using the fact that
both $G_1$ and $G_2$ have exactly five non equivalent one dimensional
$K$-representations, an elementary dimension counting argument for the remaining
irreducible representations shows that $KG_1\simeq KG_2\simeq 5\,K\oplus 24\,M_5(K)$.
By Corollary~\ref{thm:lattice}, the partial group algebras over $K$ of $G_1$
and $G_2$ are isomorphic.
\end{counter}
\end{section}

\begin{section}{Some examples}\label{sec:examples}
\noindent
We conclude the paper with the explicit calculation of some examples of
groupoid algebras $K\Gamma(G)$  in a few cases.

Following the procedure described in Section~\ref{sec:structure}, 
the idea for calculating $K\Gamma(G)$ is to look at the stabilizer
of each vertex at each level of the graph $E_{\Gamma(G)}$. 
If $H$ is a subgroup of $G$ and $m\le(G:H)$ is a positive integer, 
$n(H,m)$ will denote the number of vertices of $E_{\Gamma(G)}$ at level
$k=m\times\vert H\vert$ whose stabilizer is $H$. Then, the direct summand
of $K\Gamma(G)$ corresponding to these vertices is:
\begin{equation}\label{eq:ref1}
\frac{n(H,m)}mM_m(KH).
\end{equation}
Observe that, for all finite groups $G$, it is $n(\{e\},1)=n(G,1)=1$;
in particular, $K\Gamma(G)$ contains a copy of $K$ and a copy of the group
algebra $KG$ as direct summands. 

In the following, it is useful to keep in mind that the total number of vertices
at level $k$ in $E_{\Gamma(G)}$ is:
\begin{equation}\label{eq:ref2}
\binom{\vert G\vert-1}{k-1}.
\end{equation}
\begin{Zp}\label{thm:Zp} Let $p$ be a prime number and $G=\Z_p$ the cyclic group
of order $p$. In this case, $G$ has no nontrivial subgroup, hence all the vertices
of $E_{\Gamma(G)}$ at level $k<p$ have trivial stabilizer. From (\ref{eq:ref1})
and (\ref{eq:ref2}) we obtain easily:
\begin{equation}\label{eq:kpzp}
K\Gamma(\Z_p)=\bigoplus_{k=1}^{p-1}\frac1k{\binom{p-1}{k-1}}\,M_k(K)\oplus K\Z_p.
\end{equation}
\end{Zp}

\begin{ZpZp}\label{thm:ZpZp} Let $p$ be a prime number and $G=\Z_p\times\Z_p$ be
the direct product of two cyclic groups of order $p$. The group $G$ has precisely
$p+1$ non trivial subgroups, denoted by $H_1,\ldots,H_{p+1}$, with $H_i\simeq\Z_p$
for all $i$. 

If $k$ is not a multiple of $p$, then every vertex at level $k$ in $E_{\Gamma(G)}$
has trivial stabilizer. Thus, the direct summand of $K\Gamma(G)$ corresponding
to such a level is the sum of $\frac1k\binom{p^2-1}{k-1}$ copies of the matrix
algebra $M_k(K)$.

Let's fix $k=mp$ for some $m=1,2,\ldots,p-1$. For each $i$, the number $n(H_i,m)$ is
given by the number of subsets of $G$ that can be obtained as union of $m$ distinct
(right) cosets of $H_i$ including the coset $H_i$. Hence, it is easily computed:
\begin{equation}\label{eq:numcos}
n(H_i,m)=\binom{(G:H_i)-1}{m-1}=\binom{p-1}{m-1}. 
\end{equation} 
The remaining vertices at level $mp$ have trivial stabilizer, and their number is
given by:
\begin{equation}\label{eq:numtriv}
n(\{e\},mp)=\binom{\vert G\vert-1}{mp-1}-\sum_{i=1}^{p+1}n(H_i,m)=
\binom{p^2-1}{ mp-1} -(p+1)\binom{p-1}{m-1}.
\end{equation}
Finally, we obtain:
\begin{eqnarray}\label{eq:kzpzp}
K\Gamma(\Z_p\times\Z_p)&\simeq&\bigoplus_{\substack{k=1\\
p\, \vert^{\!\!\!-}\; k}}^{p^2-1}\frac1k\binom{p^2-1}{ k-1}M_k(K)
\bigoplus_{m=1}^{p-1}\frac{p+1}m\binom{p-1}{ m-1}\,M_m(K\Z_p)\oplus\nonumber\\
&&\!\!\!\!\!\!\!\!\!\!\!\!\!\!\!\!\!\!\!\!\!\!\!
\bigoplus_{m=1}^{p-1}\frac1{mp}\left[
\binom{p^2-1}{ mp-1} -(p+1)\binom{p-1}{ m-1}\right]\,M_{mp}(K)\oplus
K(\Z_{p}\times\Z_p).
\end{eqnarray}
In particular, for $p=2$, if $K$ is algebraically closed and
${\rm char}(K)\ne2$, formula
(\ref{eq:kzpzp}) gives:
\begin{equation}\label{eq:kz2z2}
K\Gamma(\Z_2\times\Z_2)\simeq 11\,K\oplus M_3(K).
\end{equation}
\end{ZpZp}

\begin{Zpq}\label{thm:Zpq}
Let $p$ and $q$ be distinct prime numbers and let $G$ be the product group $\Z_p\times
\Z_q\simeq\Z_{pq}$. The group $G$ has precisely two non trivial
subgroups, isomorphic to the cyclic groups $\Z_p$ and $\Z_q$, which are
both maximal. Arguing as in the previous example, one computes easily
the groupoid algebra $K\Gamma(\Z_p\times\Z_q)$ as:
\begin{equation}\label{eq:zpzq}
\begin{split}
K\Gamma(\Z_p\times\Z_q)&\simeq \bigoplus_{\substack{k=1\\
p\, \vert^{\!\!\!-}\, k,\,q\, \vert^{\!\!\!-}\, k}}^{pq-1}\frac1k\binom{pq-1}{ k-1}M_k(K)
\oplus\\ &
\bigoplus_{m=1}^{p-1}\frac1m\binom{p-1}{m-1}M_m(K\Z_q)
\bigoplus_{m=1}^{q-1}\frac1m\binom{q-1}{m-1}M_m(K\Z_p)\oplus\\ &
\bigoplus_{m=1}^{q-1}\frac1{mp}\left[\binom{pq-1}{mp-1}-\binom{q-1}{m-1}\right]
M_{mp}(K)\oplus\\ &
\bigoplus_{m=1}^{p-1}\frac1{mq}\left[\binom{pq-1}{mq-1}-\binom{p-1}{m-1}\right]
M_{mq}(K)\oplus K(\Z_p\times\Z_q).
\end{split}
\end{equation}
In particular, if $p=2$, $q=3$, $K$ is algebraically closed 
and ${\rm char}(K)\ne2,3$, formula
(\ref{eq:zpzq}) gives:
\begin{equation}\label{eq:z6}
K\Gamma(\Z_6)\simeq 12\,K\oplus 4\,M_2(K)\oplus 3\,M_3(K)\oplus 2\,M_4(K)\oplus
M_5(K).
\end{equation}
\end{Zpq}

\begin{Zpp}\label{thm:Zpp} Let $p$ be a prime number and $G=\Z_{p^2}$ be the
cyclic group of order $p^2$. In this case, $G$ has only one non trivial
subgroup $H\simeq\Z_p$. If $k$ is not a multiple of $p$, then all vertices
at level $k$ in $E_{\Gamma(G)}$ have trivial stabilizer. If $k=mp$ for
some $m=1,2,\ldots,p-1$, the number $n(H,mp)$ is given by the number of
subsets of $G$ that are union of $m$ distinct (right) cosets of $H$, including
the coset $H$. Hence, it is $n(H,mp)=\binom{p-1}{ m-1}$; all the other vertices
at level $mp$, whose number is $\binom{p^2-1}{ mp-1}-\binom{p-1}{ m-1}$, have
trivial stabilizer. Therefore, we have:
\begin{eqnarray}\label{eq:kzpp}
K\Gamma(\Z_{p^2})&\simeq&
\bigoplus_{\substack{k=1\\
p\not\vert\; k}}^{p^2-1}\frac1k\binom{p^2-1}{ k-1}M_k(K)
\bigoplus_{m=1}^{p-1}\frac{1}m\binom{p-1}{ m-1}\,M_m(K\Z_p)\oplus\nonumber\\
&&\bigoplus_{m=1}^{p-1}\frac1{mp}\left[
\binom{p^2-1}{ mp-1} -\binom{p-1}{ m-1}\right]\,M_{mp}(K)\oplus K(\Z_{p^2}).
\end{eqnarray}
In particular, for $p=2$, if $K$ is algebraically closed
and ${\rm char}(K)\ne2$, formula
(\ref{eq:kzpp}) gives:
\begin{equation}\label{eq:kz22}
K\Gamma(\Z_4)\simeq 7\,K\oplus M_2(K)\oplus M_3(K).
\end{equation}
We remark that, from(\ref{eq:kz2z2}) and (\ref{eq:kz22}), we have:
\[K\Gamma(\Z_2\times\Z_2)\not\simeq K\Gamma(\Z_4)\]
if $K$ is algebraically closed and ${\rm char}(K)\ne2$.
\end{Zpp}

\begin{S3}\label{thm:S3} We now consider the non abelian group $G=S_3$, the symmetric 
group on three elements. The group $S_3$ has three subgroups of order two and
one subgroup of order three. Thus, at level $k=2$ we have three vertices in $E_{\Gamma(S_3)}$ 
with stabilizer isomorphic to $\Z_2$, and the other two vertices have trivial
stabilizer. At level $k=3$ there is one vertex with stabilizer isomorphic to $\Z_3$
and the other nine vertices have trivial stabilizer. At level $k=4$ there are six vertices
with stabilizer isomorphic to $\Z_2$ and the other four vertices have trivial stabilizer.
Finally, all the five vertices at level $k=5$ have trivial stabilizer. Summarizing, we
have:
\begin{equation}\label{eq:ks3}
\begin{split}
K\Gamma(S_3)\simeq &\, K\oplus M_2(K)\oplus 3\,M_3(K)\oplus M_4(K)\oplus M_5(K)\oplus 
\\
&\oplus 3\,K\Z_2\oplus K\Z_3\oplus3\,M_2(K\Z_2)\oplus KS_3.
\end{split}
\end{equation}
If $K$ is algebraically closed and
${\rm char}(K)\ne2,3$, then $KS_3\simeq 2K\oplus M_2(K)$,
and:
\begin{equation}\label{eq:S3}
K\Gamma(S_3)\simeq12\,K\oplus 8\,M_2(K)\oplus 3\,M_3(K)\oplus
M_4(K)\oplus M_5(K).
\end{equation}
\end{S3}
It is interesting to consider weaker questions than
the isomorphism problem, namely, to determine which properties of
$G$  can be deduced from the structure of the partial group algebra $\kpar(G)$.
For instance, it is not clear how to check the commutativity of $G$.
Observe that, from (\ref{eq:z6}) and (\ref{eq:S3}) we have that the commutative
group $\Z_6$ has partial group algebra whose center has dimension equal to
$22$, while the center of the partial group algebra of the noncommutative 
group $S_3$ has bigger dimension, equal to $25$.
\end{section}


\end{document}